# Proper Selection of Obreshkov-Like Numerical Integrators Used as Numerical Differentiators for Power System Transient Simulation

Sheng Lei
Hitachi Energy
San Jose, CA, USA
sheng.lei@hitachienergy.com

Alexander Flueck
Illinois Institute of Technology
Chicago, IL, USA
flueck@iit.edu

*Abstract*—Obreshkov-like numerical integrators have been widely applied to power system transient simulation. Misuse of the numerical integrators as numerical differentiators may lead to numerical oscillation or bias. Criteria for Obreshkov-like numerical integrators to be used as numerical differentiators are proposed in this paper to avoid these misleading phenomena. The coefficients of a numerical integrator for the highest order derivative turn out to determine its suitability. Some existing Obreshkov-like numerical integrators are examined under the proposed criteria. It is revealed that the notorious numerical oscillations induced by the implicit trapezoidal method cannot always be eliminated by using the backward Euler method for a few time steps. Guided by the proposed criteria, a frequency response optimized integrator considering second order derivative is put forward which is suitable to be used as a numerical differentiator. Theoretical observations are demonstrated in time domain via case studies. The paper points out how to properly select the numerical integrators for power system transient simulation and helps to prevent their misuse.

## I. INTRODUCTION

Obreshkov (also spelled as Obrechkoff)-like numerical integrators refer to single-step or multistep ones considering higher order derivative [1], [2]. They are introduced into transient simulation of electronic circuits and power systems to simultaneously achieve efficiency and accuracy [2], [3]. In [2], numerical integrators from the original Obreshkov family [1], [4] are introduced into circuit simulation. In [5], explicit high order Taylor series method [1], [6] is introduced into power system transient stability (TS) simulation. Note that explicit high order Taylor series method can be understood as a special case of Obreshkov numerical integrators. Reference [3] generalizes the idea of modifying the frequency response of the error of the implicit trapezoidal method proposed in [7] to Obreshkov-like numerical integrators, so that new members are introduced into the Obreshkov family. Obreshkov-like numerical integrators are further applied to power system electromagnetic transient (EMT) simulation in [3]. Many commonly used numerical integrators in power system transient simulation, such as the implicit trapezoidal method and the backward Euler method, are also Obreshkov-like.

In transient simulation, numerical integrators sometimes have to be used as numerical differentiators [8]-[11]. In these situations, unpleasant or even misleading numerical phenomena, such as numerical oscillation, may be induced by certain numerical integrators. For example, the implicit trapezoidal method is notorious for inducing numerical oscillations when used as a numerical differentiator [8], [9], [11]. On the other hand, some other numerical integrators are found to never induce numerical oscillations, such as the backward Euler method and the 2-step backward differentiation formula [8], [9]. Traditionally it is believed that the numerical oscillations of the implicit trapezoidal method can be eliminated by temporarily using the backward Euler method for a few time steps [9], [12], [13].

Considering the potential hazards they may bring, it is critical to explore the suitability of Obreshkov-like numerical integrators to be used as numerical differentiators. Schemes to mitigate the hazards, if there are any, may be further designed accordingly. However, to the authors' knowledge, this topic has not yet been reported in the literature. This paper aims at investigating how to properly select Obreshkov-like numerical integrators which are suitable to be used as numerical differentiators. Main contributions of the paper are threefold. First, criteria for the numerical integrators regarding the suitability are proposed. The suitability turns out to be an intrinsic property of the numerical integrators determined by the coefficients for the highest order derivative. A sufficient and necessary condition for an ideal situation is identified where any unpleasant numerical phenomena will be rapidly eliminated. Second, the suitability of some existing numerical integrators is examined according to the proposed criteria. It is revealed that the numerical oscillations induced by the implicit trapezoidal method cannot always be eliminated by temporarily using the backward Euler method for a few time steps. Third, a frequency response optimized integrator considering second order derivative which is suitable to be used as a numerical differentiator is put forward. Theoretical observations are demonstrated from the time domain perspective with case studies.

The rest of the paper is organized as follows. Section II proposes the criteria for Obreshkov-like numerical integrators to be used as numerical differentiators. In Section III, suitability of some existing numerical integrators is examined while the novel numerical integrator mentioned above is put forward. Time domain case studies demonstrating the theoretical observations are presented in Section IV. Section V concludes the paper.

## II. CRITERIA FOR OBRESHKOV-LIKE NUMERICAL INTEGRATORS TO BE USED AS NUMERICAL DIFFERENTIATORS

Suppose that $u$ is a time-dependent input. $u$ and its derivatives up to the $(k\text{-}1)^{\text{th}}$ order are given by discrete points in time with a uniform time step $h$. More specifically, the following information is given

$$\begin{cases} u_t, u_{t-h}, u_{t-2h}, u_{t-3h} \cdots \\ u'_t, u'_{t-h}, u'_{t-2h}, u'_{t-3h} \cdots \\ \cdots \\ u^{(k-1)}_t, u^{(k-1)}_{t-h}, u^{(k-1)}_{t-2h}, u^{(k-1)}_{t-3h} \cdots \end{cases} \quad (1)$$

where $t$ is the current time instant; $h$ is a positive real number; $k$ is a positive integer. The objective is to calculate $u^{(k)}_t$, a numerical approximation to the $k^{\text{th}}$ order derivative of $u$ at the current time instant. To this end, a numerical differentiator [6] may be used.

The collection of $u^{(k)}_t$, $u^{(k)}_{t\text{-}h}$, $u^{(k)}_{t\text{-}2h}$, ... forms a discrete-time variable $\bar{u}^{(k)}$, which is a numerical approximation to $u^{(k)}$. If the values of $\bar{u}^{(k)}$ at previous time steps are to be utilized in the calculation of $u^{(k)}_t$, which become available at the current time instant, an Obreshkov-like numerical integrator can be adopted. In particular, applying an Obreshkov-like numerical integrator, $u_t$ is expressed as

$$u_t = \sum_{j=1}^{m} c^0_{-j} u_{t-jh} + \sum_{i=1}^{k} \sum_{j=0}^{m} c^i_{-j} u^{(i)}_{t-jh} \quad (2)$$

where $c^i_{-j}$ are coefficients to be determined. A specific selection for these coefficients determines a numerical integrator. $u^{(i)}$ denotes the $i^{\text{th}}$ order derivative of $u$. $m$ is the maximum number of previous time steps considered, which is a positive integer. (2) can be converted into

$$u^{(k)}_t = -\sum_{j=1}^{m} \frac{c^k_{-j}}{c^k_0} u^{(k)}_{t-jh} + \frac{1}{c^k_0} u_t - \sum_{j=1}^{m} \frac{c^0_{-j}}{c^k_0} u_{t-jh} - \sum_{i=1}^{k-1} \sum_{j=0}^{m} \frac{c^i_{-j}}{c^k_0} u^{(i)}_{t-jh} \quad (3)$$

Note that $c^k_0$ has to be nonzero. In this situation, the numerical integrator is used as a numerical differentiator.

Further examination on (3) will reveal the potential hazards that may be induced by using a numerical integrator as a numerical differentiator. Introducing the auxiliary variables

$$y_{j,t} = u^{(k)}_{t-jh} \quad (j=1,\cdots,m-1) \quad (4)$$

Note that

$$y_{j,t} = y_{j-1,t-h} \quad (j=2,\cdots,m-1) \quad (5)$$

The following discrete linear system is obtained

$$\begin{pmatrix} u^{(k)} \\ y_1 \\ y_2 \\ \vdots \\ y_{m-2} \\ y_{m-1} \end{pmatrix}_t = \begin{pmatrix} -\frac{c^k_{-1}}{c^k_0} & -\frac{c^k_{-2}}{c^k_0} & -\frac{c^k_{-3}}{c^k_0} & \cdots & -\frac{c^k_{-(m-1)}}{c^k_0} & -\frac{c^k_{-m}}{c^k_0} \\ 1 & 0 & 0 & \cdots & 0 & 0 \\ 0 & 1 & 0 & \cdots & 0 & 0 \\ \vdots & \vdots & \vdots & \ddots & \vdots & \vdots \\ 0 & 0 & 0 & \cdots & 0 & 0 \\ 0 & 0 & 0 & \cdots & 1 & 0 \end{pmatrix} \begin{pmatrix} u^{(k)} \\ y_1 \\ y_2 \\ \vdots \\ y_{m-2} \\ y_{m-1} \end{pmatrix}_{t-h} + EIU \quad (6)$$

where $EIU$ denotes an expression involving $u$ and its derivatives up to the $(k\text{-}1)^{\text{th}}$ order, the details of which are not of interest.

Equation (6) defines a recursive calculating process. The initial value of $\bar{u}^{(k)}$ has to be externally provided. If a multistep numerical integrator is adopted, the initial value is required at multiple time steps. If the initial value is properly given in the sense that it matches $u^{(k)}$ at the beginning at the corresponding time instants, (6) gives a reasonable numerical approximation to $u^{(k)}$ on the studied time interval. If on the other hand, the initial value is improper which does not match $u^{(k)}$ at the beginning, it can be decomposed into a proper component and an additional error. According to the superposition property of linear systems [14], the error will cause additional phenomena to the approximation to $u^{(k)}$. In fact, the notorious numerical oscillation induced by the implicit trapezoidal method used as a numerical differentiator [8], [9], [11] is of this type.

Improper provision of the initial value for $\bar{u}^{(k)}$ is in fact far from being rare in transient simulation. Discontinuity events such as a switching operation in electronic circuits and a fault in power systems will result in instantaneous change in the value of some variables including some $u^{(k)}$ [15], [16]. For example, the terminal voltage as the scaled derivative of the current going through a linear time-invariant inductor can change instantaneously. Nevertheless, the computed values of $\bar{u}^{(k)}$ are stored and not modified to reflect the instantaneous change. Instead, they are used as the improper initial value for later time steps. Having this situation in mind, the frequent encounter of unpleasant phenomena like numerical oscillations is not hard to understand.

How the error impacts the numerical results depends on the state transition matrix in (6). According to linear system theory [14], if all the eigenvalues of the matrix have a magnitude smaller than 1, the impact of the error will vanish asymptotically; otherwise the error will have permanent impact on the numerical results. Especially, if -1 is an eigenvalue, numerical oscillation will appear as a result; if 1 is an eigenvalue, bias will be reflected in transient simulation. An ideal situation is that all the eigenvalues of the state transition matrix are 0. In this case, the impact of the error will be rapidly eliminated.

Note further that the eigenvalues of the state transition matrix in (6) are exactly the roots of the polynomial

$$\lambda^m + \frac{c^k_{-1}}{c^k_0} \lambda^{m-1} + \frac{c^k_{-2}}{c^k_0} \lambda^{m-2} + \cdots + \frac{c^k_{-m}}{c^k_0} \quad (7)$$

where the coefficients are those of the numerical integrator for the highest order derivative at the time steps considered. The ideal situation where all the roots are 0 happens if and only if the coefficients satisfy

$$c_{-j}^k = 0 \ (j \neq 0) \quad (8)$$

In summary, whether or not an Obreshkov-like numerical integrator is suitable to be used as a numerical differentiator depends on the polynomial formed by its coefficients for the highest order derivative as (7). If the magnitude of all the roots is smaller than 1, the unpleasant impact from improper initial value will vanish asymptotically; otherwise the improper initial value will distort the numerical results in a permanent manner. A root at -1 will induce numerical oscillation; a root at 1 will induce bias. The ideal situation where the impact of improper initial value is rapidly eliminated happens if and only if the coefficients for the highest order derivative at the previous time steps are all 0; in this situation all the roots are 0.

It is true that undesirable numerical phenomena can be avoided if proper initial value is provided. Unfortunately this is generally unachievable due to numerical error in calculation.

III. EXAMINING SOME OBRESHKOV-LIKE NUMERICAL INTEGRATORS USED AS NUMERICAL DIFFERENTIATORS

A. Backward Euler Method and 2-Step Backward Differentiation Formula

The backward Euler method (BE) is a single-step numerical integrator considering first order derivative [1], [6]. By BE,

$$u_t = u_{t-h} + hu_t' \quad (9)$$

The 2-step backward differentiation Formula (BDF2) is a 2-step numerical integrator considering first order derivative [1], [6]. By BDF2,

$$u_t = \frac{4}{3}u_{t-h} - \frac{1}{3}u_{t-2h} + \frac{2}{3}hu_t' \quad (10)$$

From (9) and (10), it is observed that these two numerical integrators are suitable to be used as numerical differentiators because the coefficients for the highest order derivative at previous time steps are all 0. The impact of improper initial value will be dramatically eliminated.

B. Implicit Trapezoidal Method

The classical implicit trapezoidal method (TR) is a single-step numerical integrator considering first order derivative [1], [6]. By TR,

$$u_t = u_{t-h} + \frac{h}{2}u_t' + \frac{h}{2}u_{t-h}' \quad (11)$$

The polynomial (7) for TR is

$$\lambda + 1 \quad (12)$$

The root of (12) is -1. Therefore TR is not suitable to be used as a numerical differentiator. Improper initial value will cause numerical oscillations.

Traditionally, it is believed that the numerical oscillations induced by TR can be eliminated by temporarily using BE for several time steps or half time steps [9], [12], [13]. Numerical examples with a step input demonstrate this plausible argument [8], [9]. In fact, the derivative of a step input is a Dirac delta function [8], [9]; the value at the switching time instant is infinite while it is 0 elsewhere. Note that infinity cannot be strictly represented in digital simulation. Therefore it is impossible to provide a correct value at the switching time instant which serves as the proper initial value for later time steps. When using BE to numerically differentiate a step input, the calculated derivative turns to 0 in 2 time steps after the switching time instant [8], [9] which coincides with the analytical derivative. In this special case, BE happens to provide the proper initial value for TR, which triggers no numerical oscillations later. Unfortunately such coincidence cannot be generally guaranteed.

In the general setting, the calculated derivative from BE typically contains numerical error, making it different from the analytical derivative. The inaccurate calculated value will become the improper initial value for the later calculation carried out with TR, leading to numerical oscillations. In the next section, a counterexample will be presented, where the numerical oscillations induced by TR cannot be eliminated by temporary use of BE.

C. Single-Step Numerical Integrators Considering Second Order Derivative

By a single-step numerical integrator considering second order derivative,

$$u_t = c_{-1}^0 u_{t-h} + c_0^1 u_t' + c_{-1}^1 u_{t-h}' + c_0^2 u_t'' + c_{-1}^2 u_{t-h}'' \quad (13)$$

where the coefficients are to be determined.

Reference [7] proposes modifying the implicit trapezoidal method so that the frequency response of the error is set to zero at a specified nonzero angular frequency $\omega_{select}$. Such a modification makes the resulting numerical integrator especially accurate for those differential equations of which the state variable is dominated by the frequency component around $\omega_{select}$. This idea can be extended to (13). Specifically, performing the Laplace transform on both sides of (13)

$$U = c_{-1}^0 U e^{-sh} + c_0^1 sU + c_{-1}^1 sUe^{-sh} + c_0^2 s^2 U + c_{-1}^2 s^2 U e^{-sh} \quad (14)$$

The $s$-domain error expression is

$$U - (c_{-1}^0 U e^{-sh} + c_0^1 sU + c_{-1}^1 sUe^{-sh} + c_0^2 s^2 U + c_{-1}^2 s^2 U e^{-sh}) \quad (15)$$

The $s$-domain relative error expression is

$$1 - (c_{-1}^0 e^{-sh} + c_0^1 s + c_{-1}^1 se^{-sh} + c_0^2 s^2 + c_{-1}^2 s^2 e^{-sh}) \quad (16)$$

If the coefficients are chosen so that $j\omega_{select}$ and $-j\omega_{select}$ are a pair of roots of (16), then the frequency response of the error is set to zero at the angular frequency $\omega_{select}$, making the numerical integrator accurate for signals with a dominant frequency component at $\omega_{select}$. If 0 is a root of (16), then the frequency response of the error is set to zero at 0 Hz, leading to high accuracy for slow variants. Note that a multiple root introduces smaller error around the specified frequency. Desirable selection is achieved by solving equations regarding these root conditions.

According to the proposed criteria for Obreshkov-like numerical integrators used as numerical differentiators, the coefficient for the highest order derivative at the previous time step $c_{-1}^2$ should be 0 to rapidly eliminate the impact of any improper initial value.

Considering the above mentioned conditions, this paper proposes a numerical integrator as follows which is suitable to be used as a numerical differentiator

$$c_{-1}^0 = 1, \quad c_{-1}^2 = 0$$
$$c_0^1 = \frac{-(\sin(\omega_{select}h) - \omega_{select}h\cos(\omega_{select}h))}{\omega_{select}(\cos(\omega_{select}h) - 1)}$$
$$c_{-1}^1 = \frac{\sin(\omega_{select}h) - \omega_{select}h}{\omega_{select}(\cos(\omega_{select}h) - 1)} \quad (17)$$
$$c_0^2 = \frac{-(2\cos(\omega_{select}h) + \omega_{select}h\sin(\omega_{select}h) - 2)}{\omega_{select}^2(\cos(\omega_{select}h) - 1)}$$

With this set of coefficients, $j\omega_{select}$ and $-j\omega_{select}$ are a single root and 0 is a double root of the relative error expression (16). This numerical integrator is referred to as Integrator E in this paper.

Other single-step numerical integrators considering second order derivative to be examined in this paper are listed in Table I. Integrators A and B are defined in [3]. Integrator C can be found in [1] and [4]. Integrator D can be found in [6]. Integrator F can be found in [4]. Interested readers may refer to the references for more details. In Table I, the coefficient $c_0^2$ of Integrator A is

$$-\frac{1}{\omega_{select}^2} + \frac{h}{2\omega_{select}}\cot(\frac{\omega_{select}h}{2}) \quad (18)$$

The numerical integrators are examined in Table II to see whether they are suitable to be used as numerical differentiators.

TABLE I. OTHER NUMERICAL INTEGRATORS CONSIDERING SECOND ORDER DERIVATIVE TO BE EXAMINED IN THIS PAPER

| Coefficient | Numerical Integrator | | | | |
|---|---|---|---|---|---|
| | A | B | C | D | F |
| $c_{-1}^0$ | 1 | 1 | 1 | 1 | 1 |
| $c_0^1$ | $\frac{h}{2}$ | $\frac{sin(\omega_{select}h)}{\omega_{select}}$ | $\frac{h}{2}$ | $h$ | $\frac{2}{3}h$ |
| $c_{-1}^1$ | $\frac{h}{2}$ | 0 | $\frac{h}{2}$ | 0 | $\frac{1}{3}h$ |
| $c_0^2$ | (18) | $\frac{cos(\omega_{select}h)-1}{\omega_{select}^2}$ | $-\frac{h^2}{12}$ | $-\frac{h^2}{2}$ | $-\frac{1}{6}h^2$ |
| $c_{-1}^2$ | $-c_0^2$ | 0 | $\frac{h^2}{12}$ | 0 | 0 |
| Roots at $\pm j\omega_{select}$ of (16) | 1 | 1 | 0 | 0 | 0 |
| Roots at 0 of (16) | 3 | 1 | 5 | 3 | 4 |

TABLE II. EXAMINATION OF THE NUMERICAL INTEGRATORS

| Feature | Numerical Integrator | | | | |
|---|---|---|---|---|---|
| | A | B | C | D | F |
| Polynomial (7) | $\lambda-1$ | $\lambda$ | $\lambda-1$ | $\lambda$ | $\lambda$ |
| Roots of the polynomial | 1 | 0 | 1 | 0 | 0 |
| Suitable to be used as a numerical differentiator | No | Yes | No | Yes | Yes |
| Hazards triggered by improper initial values | Bias | -- | Bias | -- | -- |

## IV. TIME DOMAIN DEMONSTRATION

In this section, it is assumed that the input $u$ is given by
$$u = \cos(\omega_{syn}t) \quad (19)$$
where $\omega_{syn} = 120\pi$ s$^{-1}$. Note that

$$u' = -\omega_{syn}\sin(\omega_{syn}t) \quad (20)$$
$$u'' = -\omega_{syn}^2\cos(\omega_{syn}t) \quad (21)$$

Theoretical observations in Section III are to be demonstrated in this section via time domain case studies, which start from 0 s. The proper initial value for $\bar{u}'$ is 0 while that for $\bar{u}''$ is $-\omega_{syn}^2$. Numerical results from numerical integrators used as numerical differentiators are compared to the analytical result. For those frequency response optimized integrators, $\omega_{select}$ is set to $\omega_{syn}$ when applicable.

### A. Failure of BE in Eliminating the Numerical Oscillations Induced by TR

In this subsection, $\bar{u}'$ is to be numerically calculated. The initial value is deliberately set to 300 to create an improper one. The numerical result from TR is compared to the analytical result in Fig. 1. A step size of 1 ms is used. As theoretically predicted in Section III, the notorious numerical oscillation is observed. In Fig. 2, the numerical results from TR with temporary use of BE at the beginning are compared to the analytical result. The same 1 ms step size is used. Two schemes are considered. The first one implements 2 half time steps of BE at the beginning; the other 4 half time steps. Both schemes fail in eliminating the numerical oscillation, as discussed in Section III.

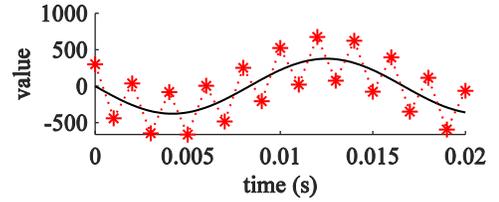

Figure 1. Numerical result from TR compared to the analytical result. Solid line: analytical result. Dotted line with asterisk markers: numerical result from TR.

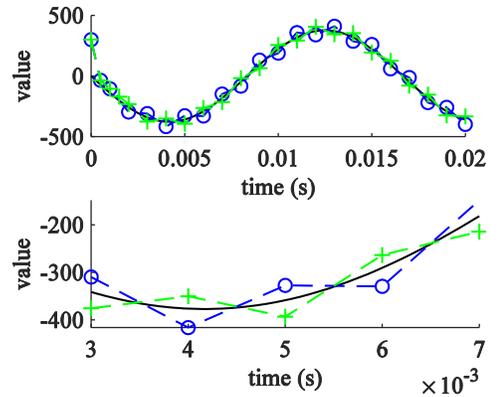

Figure 2. Numerical results from TR with temporary use of BE at the beginning compared to the analytical result. Solid line: analytical result. Dashed line with circle markers: 2 half time steps of BE are used. Dashed line with plus sign markers: 4 half time steps of BE are used. Top: the comparison from 0 to 0.02 s. Bottom: zooming-in at 0.005 s.

### B. Numerical Integrators Considering Second Order Derivative

In this subsection, $\bar{u}''$ is to be calculated by Integrators A-F. The initial value is deliberately set to 0 to create an improper one. In Fig. 3, numerical results from Integrators A, C and E are compared to the analytical result. A step size of 2 ms is used. For better visualization, other numerical integrators are not included. Fig. 3 demonstrates the theoretical predictions in Section III. Integrators A and C induce significant bias. On the other hand, the proposed Integrator E eliminates the impact of the improper initial value rapidly.

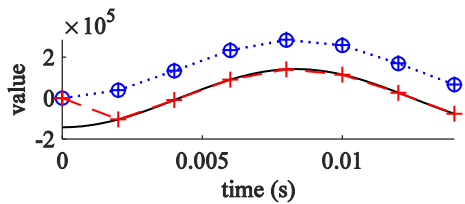

Figure 3. Numerical results from numerical integrators considering second order derivative compared to the analytical result. Solid line: analytical result. Dotted line with circle markers: result from Integrator A. Dotted line with plus sign markers: result from Integrator C. Dashed line with plus sign markers: result from Integrator E.

In order to quantitatively study the numerical error of each numerical integrator, the following error measurement is used. The relative error regarding $\bar{u}''$ from a numerical integrator used as a numerical differentiator with a specified step size is defined as

$$err(\bar{u}'') = \|\bar{u}'' - u''\|_2 / \|u''\|_2 \times 100 \quad (22)$$

Note that numerical results are discrete; the 2-norm is thus calculated at common time instants of $\bar{u}''$ and $u''$. Table III lists the error of Integrators B, D, E and F with different step sizes. Integrators A and C are not included because they induce bias and are not suitable to be used as numerical differentiators. Computations performed to generate Table III are from 0 to 1 s. The first two time steps are not considered to exclude the impact of improper initial value.

TABLE III. Error of Numerical Integrators Used as Numerical Differentiators with Different Step Sizes

| Step Size (μs) | Numerical Integrator | | | |
|---|---|---|---|---|
| | B | D | E | F |
| 125 | 0.0000 | 1.5709 | 0.0000 | 0.0185 |
| 250 | 0.0000 | 3.1418 | 0.0000 | 0.0740 |
| 500 | 0.0000 | 6.2820 | 0.0000 | 0.2959 |
| 1000 | 0.0000 | 12.5428 | 0.0000 | 1.1809 |
| 2000 | 0.0000 | 24.8785 | 0.0000 | 4.6812 |
| 4000 | 0.0000 | 48.0113 | 0.0000 | 18.0758 |

From Table III it is learnt that Integrators B and E are exactly accurate despite the step size. Their accuracy in time domain matches the theoretical prediction in frequency domain presented in Section III-C. According to the frequency domain error analysis, Integrator E is in fact more accurate than Integrator B in general. Integrator D is the least accurate; if the step size is doubled, its error is also roughly doubled. For Integrator F, if the step size is doubled, the error increases by a factor of around 4. Integrator F is much more accurate than Integrator D in time domain, which can also be told by their frequency response of the error. The proposed Integrator E is very promising as a numerical differentiator for sinusoidal inputs considering its suitability and high accuracy.

## V. Conclusion

Undesirable numerical phenomena induced by numerical integrators used as numerical differentiators are attributed to the intrinsic suitability of the numerical integrators and triggered by improper initial value. The coefficients for the highest order derivative determine the suitability. The common practice of using the backward Euler method for a few time steps in an attempt to eliminate the numerical oscillations induced by the implicit trapezoidal method is unreliable. A frequency response optimized integrator considering second order derivative is proposed which is suitable to be used as a numerical differentiator.